\documentclass[11pt]{amsart}
\usepackage[leqno]{amsmath}
\usepackage{amsfonts, amssymb,srcltx, amsopn, color}

\usepackage{graphicx}
\usepackage{xspace}

\usepackage{epsfig}
\usepackage{color}
\usepackage{amsmath}
\usepackage{amssymb}

\newtheorem{theorem}{Theorem}
\newtheorem{corollary}{Corollary}

\theoremstyle{definition}

\usepackage{color}
\newcommand{\commentout}[1]{}

\usepackage{graphicx}

\begin{document}

\thispagestyle{empty}

\centerline{\Large\bf Distance--preserving subgraphs of Johnson graphs}

\vspace{6mm}

\centerline{\sc Victor Chepoi}

\vspace{3mm}

\medskip
\begin{small}
\medskip
\centerline{Laboratoire d'Informatique Fondamentale, Aix-Marseille Universit\'e and CNRS,}
\centerline{Facult\'e des Sciences de Luminy, F-13288 Marseille Cedex 9, France}
\centerline{\texttt{victor.chepoi@lif.univ-mrs.fr}}
\end{small}

\bigskip\bigskip\noindent{\bf Abstract.}  {\footnotesize  We give a characterization of distance--preserving subgraphs of
Johnson graphs, i.e. of graphs which are isometrically embeddable into Johnson graphs (the Johnson graph $J(m,\Lambda)$
has the subsets of cardinality $m$ of a set $\Lambda$ as the vertex--set and two such sets $A,B$ are adjacent iff
$|A\triangle B|=2$). Our characterization is similar
to the characterization of D. \v{Z}. Djokovi\'c (J. Combin. Th. Ser. B  14 (1973), 263--267) of distance--preserving subgraphs of hypercubes
and provides an explicit description of the wallspace (split system) defining the embedding.}

\thispagestyle{empty}

\section{Introduction}\label{s:intro}

\subsection{Basic notions} All graphs considered in this note are undirected, connected, contain
no multiple edges or loops, but are not necessarily finite.
The {\it distance} $d_G(u,v)$ between two vertices $u$ and $v$ of a graph $G=(V,E)$
is the length of a shortest $(u,v)$--path. The {\it interval}
$I(u,v)$ between $u$ and $v$ consists of all vertices on shortest
$(u,v)$--paths, that is, of all vertices (metrically) {\it between} $u$
and $v$: $I(u,v)=\{ x\in V: d(u,x)+d(x,v)=d(u,v)\}.$ A subset of vertices $A$ of $G$
is called {\it  convex} if it includes the interval between any pair of its
vertices.  A graph $G=(V,E)$
is {\it isometrically embeddable} into a graph $G'=(V',E')$
if there exists a mapping $\varphi : V\rightarrow V'$ such that $d_{G'}(\varphi (u),\varphi
(v))=d_G(u,v)$ for all vertices $u,v\in V$. The image of $G$ under $\varphi$ is called an
{\it isometric} or a {\it distance--preserving subgraph} of $G'$. More generally, for a positive
integer $k$,  a {\it scale $k$ embedding}
of a graph $G$ into a graph $G'$ is a mapping $\varphi : V\rightarrow V'$ such that $d_{G'}(\varphi (u),\varphi
(v))=k\cdot d_G(u,v)$ for all  $u,v\in V$.

A {\it hypercube} $H(\Lambda)$ is a graph having the finite subsets of a set  $\Lambda$ as vertices
and two  such sets $A,B$ are adjacent in $H(\Lambda)$
iff $|A\triangle B|=1$. A {\it half--cube} $\frac{1}{2}H(\Lambda)$ has the finite subsets of
$\Lambda$ of even cardinality as vertices and two  such vertices $A,B$ are adjacent in $\frac{1}{2}H(\Lambda)$
iff $|A\triangle B|=2$ (analogously is defined the half--cube $\frac{1}{2}H(\Lambda)$ for finite subsets of odd size).
For an integer $m>0$, the {\it Johnson graph}  $J(m,\Lambda)$ has the
subsets of $\Lambda$ of size $m$ as vertices and
two such vertices $A,B$ are adjacent iff $|A\triangle B|=2$. Obviously, all Johnson graphs $J(m,\Lambda)$ are
isometric subgraphs of the corresponding half--cube $\frac{1}{2}H(\Lambda)$. Notice also that the half--cube
$\frac{1}{2}H(\Lambda)$ and the Johnson graphs $J(m,\Lambda)$ are scale 2 embedded in the hypercube $H(\Lambda)$.

\subsection{Distance--preserving subgraphs of hypercubes} Djokovi\'{c} \cite{Djokovic}   characterized distance--preserving
subgraphs of hypercubes in the following simple but pretty way:
{\it a graph $G=(V,E)$ can be isometrically embedded into a hypercube iff $G$ is bipartite and for any edge
$uv$, the (disjoint) sets $W(u,v)$ and $W(v,u)$, where
\begin{align*}
W(u,v)=\{ x\in V: d_G(x,u)<d_G(x,v)\} \tag{1}
\end{align*}
are convex.} If $G$ is bipartite, then $W(u,v)\cup W(v,u)=V$, whence $W(u,v)$ and $W(v,u)$ are complementary convex subsets
of $G$, called {\it convex halfspaces}. To establish an isometric embedding of $G$ into a hypercube,  Djokovi\'{c} \cite{Djokovic}
introduces the following binary relation $\theta$ on the edges of $G$:  for two edges $e=uv$ and $e'=u'v'$ we set $e\theta e'$ iff
$u'\in W(u,v)$ and $v'\in W(v,u)$. Under the conditions of the theorem, it can be shown that $e\theta e'$ iff
$W(u,v)=W(u',v')$ and $W(v,u)=W(v',u')$, whence $\theta$ is an equivalence relation. Let ${\mathcal E}=\{ E_i: i\in \Lambda\}$ be the equivalence classes
of $\theta$ and let $b$ be an arbitrary fixed vertex taken as the basepoint of $G$. For an equivalence class $E_i\in {\mathcal E}$,
let ${\mathcal W}_i=\{ H^-_i,H^+_i\}$ be the pair of complementary
convex halfspaces  of $G$ defined by setting $H^-_i:=W(u,v)$ and $H^+_i:=W(v,u)$ for an arbitrary edge $uv\in E_i$ with $b\in W(u,v)$.
Then the isometric embedding $\varphi$ of $G$ into the hypercube $H(\Lambda)$ is obtained by setting  $\varphi(v):=\{ i\in \Lambda: v\in H^+_i\}$
for any vertex $v\in V$ and ${\frak W}=\{ {\mathcal W}_i: i\in \Lambda\}$ is the signed wallspace defining this embedding (see below).

In nonbipartite graphs, for an edge $uv$ the sets $W(u,v)$ and
$W(v,u)$ no longer partition $V(G)$, therefore the set
\begin{align*}
W_=(u,v):=\{ x\in V(G): d_G(x,u)=d_G(x,v)\} \tag{2}
\end{align*}
can be nonempty. Answering a question
of Winkler \cite{Winkler_product}, the following Djokovi\'c--style characterization of
distance--preserving subgraphs of Hamming graphs (Cartesian products of complete graphs) was provided in \cite{Ch_Hamming}: {\it a graph $G$ is
isometrically embeddable into a Hamming graph iff for any edge $uv$ of $G$ the sets $W(u,v),W(v,u),W_=(u,v),$ and
their complements are convex} (other characterizations were obtained in \cite{Wil} and  \cite{Ch_Hamming}).

\subsection{Djokovi\'{c}--Winkler relation and canonical metric representation} Winkler \cite{Winkler_product} extended the Djokovi\'c relation $\theta$ to all (not necessarily
bipartite) graphs $G$: $e\theta e'$ for
two edges $e=uv$ and $e'=u'v'$ iff $d(u,u')+d(v,v')\ne d(u,v')+d(v,u')$. 
He proved that {\it $\theta$ is transitive iff $G$ isometrically
embeds into a Cartesian product of $K_3$'s.}    Djokovi\'c--Winkler relation
$\theta$ plays a determinant role in the Graham and Winkler's \cite{GrWi} {\it canonical metric representation of graphs}: it was shown in \cite{GrWi}
that {\it any  finite connected graph $G$ has a unique isometric embedding $\varphi$ into the Cartesian product $\Pi_{i=1}^m G_{i}$ in which each factor $G_{i}$ is
irreducible (i.e., not further decomposable this way)}. The factors and the embedding $\varphi$ are defined in the following way. Let  ${\mathcal E}=\{ E_i: i\in \Lambda\}$ be
the equivalence classes of the transitive closure $\theta^*$ of the Djokovi\'c--Winkler relation $\theta$ and let $m:=|\Lambda|$. For each equivalence class
$E_{i}\in {\mathcal E}$, let $G_{i}$ be the graph whose vertex--set is the set of connected components of the graph $G^*_{i}=(V, E\setminus E_{i})$
($G^*_{i}$ is the graph obtained from $G$ by removing the edges of $E_{i}$) and two such components are adjacent in $G_{i}$ if at least one edge of
$E_{i}$ has its ends in both components. For a vertex $v$ of $G$, the $i$th coordinate of $\varphi(v)$ is the connected component of $G^*_i$ containing $v$.

\subsection{Shpectorov characterization of $\ell_1$--graphs}
The Djokovi\'c--Winkler relation $\theta$ also played  a significant role in the beautiful Shpectorov's  proof of the  characterization of $\ell_1$--graphs
\cite{Shp} (see also \cite[Chapter 21]{DeLa}). $\ell_1$--{\it Graphs} are the graphs which can be isometrically
embedded into an $\ell_1$--space; equivalently, finite $\ell_1$--graphs are the graphs which admit a scale embedding into a hypercube \cite[Proposition 4.3.8]{DeLa}.
Shpectorov \cite{Shp} proved that {\it a finite graph $G$ is an $\ell_1$--graph iff $G$ isometrically embeds in a Cartesian product of hyperoctahedra
(complete graphs $K_{2m}$  minus a perfect matching) and half--cubes.} Equivalently, the $\ell_1$--graphs are exactly the graphs $G$ for which
the irreducible factors in
the Graham-Winkler canonical representation are either induced subgraphs of hyperoctahedra or isometric subgraphs of half--cubes. While the subgraphs of hyperoctahedra
can be easily characterized, a structural characterization of distance--preserving subgraphs of half--cubes is still missing, see \cite[Problem 21.4.1]{DeLa}.
Nevertheless, it can be tested in polynomial time if a graph $G$ is isometrically embeddable in a half--cube \cite{Shp,DeSh}. Notice also that Shpectorov's characterization 
of $\ell_1$-graphs can be viewed as a sharpenning of the characterization in the same vein of hypermetric graphs obtained by Tervilliger and Deza \cite{TeDe}. 

\subsection{Atom graph} An important ingredient in the proof of Shpectorov's theorem and in the recognition algorithm of \cite{DeSh} is the notion of the
atom graph \cite{Shp}, which can be defined as follows. Given two edges $e=uv$ and $e'=u'v'$ of $G$, set
\begin{align*}
<e,e'>:=d_G(u,v')+d_G(v,u')-d_G(u,u')-d_G(v,v').\tag{3}
\end{align*}
The quantity $<e,e'>$ takes the values $0,\pm 1,\pm 2$. Set $e\theta_1 e'$ iff $<e,e'>=\pm 2$ and set $e\theta_2 e'$ iff $<e,e'>=\pm 1.$
Note that $e\theta_1 e'$ coincides with the original  Djokovi\'c relation in the sense that $e\theta_1 e'$ iff $u'\in W(u,v)$ and $v'\in W(v,u)$.
Let $b$ be an arbitrary vertex of $G$, which we consider as a {\it basepoint}.
We call an edge $e=uv$ of $G$ {\it vertical} (with respect to $b$) if $d_G(b,u)<d_G(b,v)$ and
denote by $E_0$ the set of all vertical edges of $G$.  In the assumption that  $\theta_1$ is an equivalence
relation on $E_0$,  the {\it atom graph} $\Sigma$ (which depends of the basepoint $b$) has the set  ${\mathcal E}:=\{ E_{i}: i\in \Upsilon\}$ 
of all equivalence classes of $\theta_1$
as the vertex--set and two classes
$E_i$ and $E_j$ are adjacent in $\Sigma$  iff $e\theta_2 e'$ for all edges $e\in E_i$ and $e'\in E_j$. If $G$ is a distance--preserving subgraph
of a half--cube, then $\Sigma$ is well--defined and is a line--graph  (recall
that the {\it line--graph} $L(D)$ of a graph $D$ has the edges of $D$ as the vertex--set and two edges $e,e'$ of $D$ are adjacent in $L(D)$
iff $e$ and $e'$ share a common end). This property is crucial in the construction of  the isometric embedding of $G$ into a half--cube; for details
see \cite{DeLa,DeSh,Shp}.

\subsection{Espaces \`a murs} Let $\varphi$ be a scale $k$ embedding of a graph $G=(V,E)$ into a hypercube $H(\Lambda)$. For  $i\in \Lambda,$ let $H^-_i=\{ v\in V: i\notin \varphi(v)\}$ and
$H^+_i=\{ v\in V: i\in \varphi(v)\}$.  Then ${\mathcal W}_i=(H^-_i,H^+_i)$ is a pair of complementary halfspaces of $G$ ($H^-_i\cap H^+_i=\emptyset$ and $H^-_i\cup H^+_i=V$),
called a {\it cut} in combinatorial optimization \cite{DeLa}, a {\it split} in phylogenetic combinatorics \cite{BaDr,DrHuKoMoSp}, and a {\it wall} in geometric
group theory \cite{HaPa}.  In this paper, we will adopt the last name. A {\it wallspace} is a set $\Omega$, together with a collection $\frak W=\{ {\mathcal W}_i=(H'_i,H''_i): i\in \Lambda\}$ of bipartitions
of $\Omega$, called {\it walls}; the two sets in each bipartition are called {\it halfspaces}. Additionally, it is required that any pair of points $x,y\in \Omega$
is separated by a finite number of walls (a wall ${\mathcal W}_i=(H'_i,H''_i)$ {\it separates}
$x$ and $y$ if $x\in H'_i$ and $y\in H''_i$).  The {\it wall--distance} $d_{\frak W}(x,y)$ between two points  $x,y\in \Omega$
is the number of walls separating $x$ and $y$. A {\it signed wallspace} is a wallspace $(\Omega, {\frak W})$ together with an orientation of its walls:
each wall ${\mathcal W}_i$ has a  {\it positive halfspace} $H^+_i$ and a {\it negative halfspace} $H^-_i$, i.e., ${\mathcal W}_i=(H^-_i,H^+_i)$. If all  halfspaces of
an wallspace $(V,\frak W)$ are convex subgraphs of $G$, then $\frak W$ is called a {\it space with convex walls}.
Then a graph $G$ is scale $k$ embeddable into a hypercube $H(\Lambda)$ iff there exists a
wallspace $(V,{\frak W})$ such that $d_{\frak W}(x,y)=k\cdot d_{G}(x,y)$ or, equivalently, if there exists a space  with convex walls
such that the ends of each edge of $G$ are separated by exactly $k$ walls \cite{BaCh_weak}. Consequently, the isometric embedding of a graph $G$
into the Johnson graph $J(m,\Lambda)$ is equivalent to the existence of a signed space with convex walls $(V,{\frak W})$ such that (a) each vertex belongs to $m$
positive halfspaces and (b) the ends of each edge  of $G$ are separated by exactly two walls.

\subsection{Main result}
We continue with the main result of this paper, which provides a Djokovi\'c--like characterization of distance--preserving subgraphs of Johnson graphs and explicitly describes 
the wallspace providing the embedding:

\begin{theorem} \label{theorem1} A graph $G$ is isometrically embeddable into a Johnson graph if and only if $G$ satisfies the following two conditions:
\begin{itemize}
\item[(WC)] {\sf (Wallspace condition)} for any edge $uv$ of $G$ the subgraph  induced by $W_=(u,v)$ contains at most two connected components $W'_=(u,v), W''_=(u,v)$ 
(which are allowed to be empty) and each of the two walls
\begin{align*}
{\mathcal W}'_{uv}:=\{ W(u,v)\cup W'_=(u,v), W(v,u)\cup W''_=(u,v)\}\\
{\mathcal W}''_{uv}:=\{ W(u,v)\cup W''_=(u,v), W(v,u)\cup W'_=(u,v)\}
\end{align*}
consists of complementary convex halfspaces of $G$;
\item[(AGC)] {\sf (Atom graph condition)}  for some (in fact, for any)  basepoint $b$ of $G$ the atom graph $\Sigma$ is the line--graph of a bipartite graph 
with at least one part of the bipartition finite.
\end{itemize}
\end{theorem}

\medskip
Basis graphs of matroids are the most important examples of distance--preserving subgraphs of Johnson graphs. A {\it matroid}
on a finite set $\Lambda$ is a collection
$\mathcal B$ of subsets of $\Lambda,$ called {\it bases},  satisfying the
following exchange property: for all $A,B\in {\mathcal B}$ and $i\in A\setminus B$
there exists $j\in B\setminus A$ such that $A\setminus
\{ i\}\cup \{ j\}\in {\mathcal B}$.  All the bases of a matroid have the same
cardinality. The {\it basis graph} of a matroid $\mathcal B$ is the
graph whose vertices are the bases of $\mathcal B$ and edges are the
pairs  $A,B$ of bases such that $\vert
A\triangle B\vert=2.$   Maurer \cite{Mau} characterized the basis graphs of matroids as {\it the graphs satisfying the
following three conditions: the interval condition (IC), the positioning condition (PC), and the link
condition (LC)}. The {\it link condition} (LC) is a local version of (AGC): it asserts that the neighbors  of
some (in fact, of any) vertex $b$ induce a line--graph of a bipartite graph. The {\it interval condition} (IC)
asserts that for any two vertices $u,v$ at distance 2,  $I(u,v)$ induces a square, a pyramid, or a
3-dimensional hyperoctahedron. Finally, the {\it positioning condition} (PC) asserts that for any basepoint $b$  and any square $u_1u_2u_3u_4$ of $G$ the equality
$d_G(b,u_1)+d_G(b,u_3)=d_G(b,u_2)+d_G(b,u_4)$ holds.  It was shown recently in \cite{ChChOs}
that (LC) is implied by the two other Maurer's conditions (IC) and (PC).  We will show
below that the wall condition (WC) implies the positioning condition (PC). Together with the results of \cite{Mau} and \cite{ChChOs}, this
leads to the following characterization of basis graphs of matroids:

\begin{corollary} \label{corollary1} $G$ is the basis graph of a matroid if and only if $G$ satisfies (WC)
and (IC).
\end{corollary}

\section{Proof of Theorem 1}

\subsection{Properties of graphs satisfying (WC)} We start with some properties of graphs $G$
satisfying (WC).  We will use the notations $W'_=(u,v)$ and $W''_=(u,v)$ for
the connected components of $W_=(u,v)$ even in the case when one or both of these sets are empty.

\medskip
(2.1) {\it For any edge $uv$ of $G$, the sets $W(u,v),W(v,u), W'_=(u,v),$ and $W''_=(u,v)$ are convex.}

\medskip
{\it Proof.} Each of the sets $W(u,v),W(v,u), W'_=(u,v), W''_=(u,v)$ is the intersection of two convex  sets, one halfspace
from each of the walls ${\mathcal W}'_{uv}$ and ${\mathcal W}''_{uv}$.  $\Box$

\medskip
(2.2) {\it  For two vertical edges $e=uv$ and $e'=u'v'$ of $G$ with $b\in W(u,v)\cap W(u',v')$, the following conditions \emph{(i)}-\emph{(iii)} are equivalent:

\noindent
\emph{(i)} $e\theta_1 e'$; ~~\emph{(ii)} $<e,e'>=2$; ~~\emph{(iii)} $W(u,v)=W(u',v')$ and $W(v,u)=W(v',u')$.

\noindent
In particular, the relation $\theta_1$ is transitive.}

\medskip
{\it Proof.} From the definition of $\theta_1$ it follows that $e\theta_1 e'$  iff $u'$ and $v'$ belongs to different sets
$W(u,v)$ and $W(v,u)$.  If we suppose that $u'\in W(v,u),$ then $v'\in W(u,v)$.  Since $b\in W(u,v)$ and $u'\in I(v',b),$
we obtain a contradiction with the convexity of $W(u,v)$. Thus  $e\theta_1 e'$ is equivalent to the inclusions $u'\in W(u,v)$ and
$v'\in W(v,u)$. This is equivalent to the equalities $d_G(u,u')=d_G(v,v')$ and $d_G(u,v')=d_G(v,u')=d_G(u,u')+1$.
By (3), we conclude that $e\theta_1 e'$ iff $<e,e'>=2$, establishing the equivalence of (i) and (ii).

Now, we prove the equivalence of (i) and (iii).  Let $e\theta_1 e'$, i.e.,  $u'\in W(u,v), v'\in W(v,u)$
and $u\in W(u',v'), v\in W(v',u')$.  We assert that $W(u,v)\subseteq W(u',v')$. Pick any vertex
$x\in  W(u,v)$. If $x\in W(v',u')$, then $v'\in I(x,u')$ and since $x,u'\in W(u,v), v'\in W(v,u)$, this contradicts the convexity of $W(u,v)$. Now suppose that
$x\in W_=(u',v'),$ say $x\in W'_=(u',v')$. Since $u\in W(u',v'), v\in W(v',u'),$ and $u\in I(v,x)$, we obtain a contradiction with the convexity of $W(v',u')\cup W'_=(u',v')$.
Hence $W(u,v)\subseteq W(u',v')$. Analogously, one can show that $W(v,u)\subseteq W(v',u'), W(u',v')\subseteq W(u,v),$ and  $W(v',u')\subseteq W(v,u)$.
Consequently, $W(u,v)=W(u',v')$ and $W(v,u)=W(v',u')$. Conversely, if $W(u,v)=W(u',v')$ and $W(v,u)=W(v',u')$, since $u,u'\in W(u,v)=W(u',v')$ and
$v,v'\in W(v,u)=W(v',u'),$ thus $e\theta_1 e'$.  $\Box$

\medskip
(2.3) {\it  For two edges $e=uv$ and $e'=u'v'$ of $G$ we have $e\theta_2 e'$
iff $u'\in W(u,v)\cup W(v,u)$ and $v'\in W_=(u,v)$ or if $v'\in W(u,v)\cup W(v,u)$ and $u'\in W_=(u,v)$.
For two vertical edges $e=uv$ and $e'=u'v'$ of $G$ with $b\in W(u,v)\cap W(u',v')$, the following conditions 
are equivalent:

\noindent
\emph{(i)} $e\theta_2 e'$; ~~\emph{(ii)} $<e,e'>=1$; ~~\emph{(iii)} $u'\in W(u,v), v'\in W_=(u,v)$ or $v'\in W(v,u), u'\in W_=(u,v).$
}

\medskip
{\it Proof.} First suppose that $e\theta_2 e'$. If $u',v'\in W_=(u,v),$ then $d_G(u',u)=d_G(u',v)$ and $d_G(v',u)=d_G(v',v),$ whence $<e,e'>=0.$ Hence we can assume
that $v'\in W(v,u).$ If $u'\in W(u,v)$, then (2.2) implies $e\theta_1 e'$. Finally, if $u'\in W(v,u),$ then $d_G(u,u')=d_G(v,u')+1$ and $d_G(u,v')=d_G(v,v')+1,$ whence
$<e,e'>=0.$ Therefore, $e\theta_2e'$ and $v'\in W(v,u)$ imply  that $u'\in W_=(u,v)$. Conversely, suppose that $v'\in W(v,u)$ and $u'\in W_=(u,v)$. Then
$d_G(u,u')=d_G(v,u')$ and $d_G(u,v')=d_G(v,v')+1,$ whence $<e,e'>=1,$ yielding  $e\theta_2e'$.

Now suppose that $e=uv$ and $e'=u'v'$ are two vertical edges of $G$ with $b\in W(u,v)\cap W(u',v')$. In view of the first assertion, to establish the equivalence
of (i) and (iii), it suffices to show that if  $e\theta_2 e'$, then
$u'\notin W(v,u)$ and $v'\notin W(u,v)$. Indeed, if $u'\in W(v,u),$ then by the first assertion we infer that $v'\in W_=(u,v)$. Since $b\in W(u,v)$ and
$u'\in I(b,v')$, we obtain a contradiction with the convexity of $W(u,v)\cup W'_=(u,v),$ where $W'_=(u,v)$ is the component of $W_=(u,v)$ containing $v'$.
Analogously, if $v'\in W(u,v)$, then $u'\in W_=(u,v)$. Since $b,v'\in W(u,v)$ and $u'\in I(b,v')$ we obtain a contradiction with the convexity of $W(u,v)$.
It remains to show that $e\theta_2 e'$ implies $<e,e'>=1$ (the converse implication follows from the definition of $\theta_2$). If $u'\in W(u,v)$ and $v'\in W_=(u,v)$, then $d_G(v,u')=d_G(u,u')+1$ and $d_G(u,v')=d_G(v,v'),$
yielding $<e,e'>=1$.  Analogously, if $v'\in W(v,u)$ and $u'\in W_=(u,v)$, then $d_G(u,v')=d_G(v,v')+1, d_G(u,u')=d_G(v,u')$ and $<e,e'>=1$ again.
$\Box$

\medskip
From (2.2), $\theta_1$ is an equivalence relation on the set $E_0$ of vertical edges of $G$; let ${\mathcal E}:=\{ E_i: i\in \Lambda\}$
denote the equivalence classes of $\theta_1$ on $E_0$. Hence the atom graph $\Sigma$ introduced in Subsection 1.5 is well--defined.
In fact, we show next that the inclusion of a couple $E_i,E_j\in {\mathcal E}$ as an edge of  $\Sigma$
depends only of  the value of $<e,e'>$ of an arbitrary pair of edges $e\in E_i$ and $e'\in E_j$:

\medskip
(2.4) {\it $E_i$ and $E_j$ are adjacent in  $\Sigma$ iff  $<e,e'>=1$ for a pair of  edges $e=uv\in E_i$
and $e'=u'v'\in E_j$ with $b\in W(u,v)\cap W(u',v')$.}

\medskip
{\it Proof.} Pick any edge $e''=u''v''\in E_i$ and let $b\in W(u'',v'')$.  By (2.2), $W(u'',v'')=W(u,v)$ and $W(v'',u'')=W(v,u)$, whence $W'_=(u'',v'')=W'_=(u,v)$.
By (2.3),  either  $v'\in W(v,u)$ and $u'\in W'_=(u,v)$ or $u'\in W(u,v)$ and $v'\in W_=(u,v)$.
In the first case,  $v'\in W(v'',u'')$ and $u'\in W'_=(u'',v'')$. In the second case, $u'\in W(u'',v'')$ and $v'\in W_=(u'',v'')$. Hence $e''\theta_2e'$
by the first assertion and $<e'',e'>=1$ by the second assertion of (2.3).  By symmetry, if $e''$ is an edge of $E_j$, then $<e,e''>=1$.
Consequently, $E_i$ and $E_j$ are adjacent in $\Sigma$. $\Box$

%

\subsection{Necessity}  Next we will prove that any graph $G$ isometrically embeddable into a Johnson graph $J(m,\Lambda)$ satisfies the conditions (WC) and (AGC). Let $\varphi$ be an isometric embedding of $G$ into $J(m,\Lambda)$. Then $\varphi$ is also a scale 2 embedding of $G$ into the hypercube $H(\Lambda)$.
For a vertex $x$ of $G$, we denote by the uppercase $X$ the finite set encoding $x$, i.e., $X:=\varphi(x)$. For an edge $e=uv$ of $G$, set $\gamma(e):=U\triangle V$.
Since $|U|=|V|=m,$ there exist $i,j\in \Lambda$ such that $\gamma (e)=\{ i,j\}$ and $V=U-i+j$ (this notation stands for
$U\setminus \{ i\}\cup \{ j\}$). For a coordinate $i\in \Lambda$, let $H^-_i(\Lambda)$ and  $H^+_i(\Lambda)$ denote the set of all vertices $A$  of the hypercube $H(\Lambda)$ such that $i\notin A$ and  such that $i\in A$, respectively.

\medskip
(2.5) {\it $G^-_i:=V(G)\cap H^-_i(\Lambda)$ and $G^+_i:=V(G)\cap H^-_i(\Lambda)$ are convex halfspaces of $G$.}

\medskip
{\it Proof.} $H^-_i(\Lambda)$ and $H^+_i(\Lambda)$ are complementary convex halfspaces of $H(\Lambda)$ (as subcubes of $H(\Lambda)$).
Since $G$ is scale 2 embedded in $H(\Lambda)$,
$G^-_i$ and $G^+_i$ are complementary convex halfspaces of the graph $G$. Indeed, if say $x,y\in G^-_i$ and
$z$ is on a shortest $(x,y)$-path of $G$, then $Z$ is on a shortest $(X,Y)$-path of $H(\Lambda)$. Since $X,Y\in H^-_i(\Lambda)$, we conclude that
$Z\in H^-_i(\Lambda)$, whence $z\in G^-_i$. $\Box$

Now, let $e=uv$ be an arbitrary edge of $G$. Let  $\gamma (e)=\{ i,j\}$,  where $V=U-i+j$.

\medskip
(2.6) {\it $W(u,v)=G^+_i\cap G^-_j$ and $W(v,u)=G^-_i\cap G^+_j$. In particular, $W(u,v)$ and $W(v,u)$ are convex.}

\medskip
{\it Proof.} Pick any vertex $x\in G^+_i\cap G^-_j$. Then $X\triangle V=X\triangle U\cup \{ i,j\}$, i.e., $d_{H(\Lambda)}(X,V)=d_{H(\Lambda)}(X,U)+2.$ Consequently, $d_G(x,v)=d_G(x,u)+1,$
i.e., $x\in W(u,v)$. Hence $G^+_i\cap G^-_j\subseteq W(u,v)$. Conversely, pick any vertex $x\in W(u,v)$. Then $d_G(x,v)=d_G(x,u)+1$ and $d_{H(\Lambda)}(X,V)=d_{H(\Lambda)}(X,U)+2.$ Since
$u\in I(v,x)$ in $G$, $U$ belongs to the interval  between $V$ and $X$ of the hypercube $H(\Lambda)$.  This means that $V\cap X\subseteq U\subseteq V\cup X$. Since $i\in U\setminus V$ and $j\in V\setminus U$, this implies that $i\in X$ and $j\notin X$, whence $x\in G^+_i\cap G^-_j$. This establishes the converse inclusion $W(u,v)\subseteq G^+_i\cap G^-_j.$ Hence $W(u,v)=G^+_i\cap G^-_j$ and $W(v,u)=G^-_i\cap G^+_j$. Since $G^-_i,G^+_i, G^-_j$, and $G^+_j$ are convex, the sets $W(u,v)$ and $W(v,u)$ are convex as well. $\Box$

\medskip
(2.7) {\it Each of the sets $G^+_i\cap G^+_j$ and $G^-_i\cap G^-_j$ is  empty or induces a connected component of the subgraph $G(W_=(u,v))$ of $G$ induced by $W_=(u,v)$.}

\medskip
{\it Proof.} From (2.6) we conclude that $W_=(u,v)=(G^+_i\cap G^+_j)\cup (G^-_i\cap G^-_j)$. Each of the sets $G^+_i\cap G^+_j$ and $G^-_i\cap G^-_j$ is convex, thus
induces a connected subgraph of the graph $G$ and therefore of $G(W_=(u,v))$. Thus $G(W_=(u,v))$ contains at most
two connected components.
It remains to show that if both $G^+_i\cap G^+_j$ and $G^-_i\cap G^-_j$ are nonempty, then there is no edge between a vertex $x\in  G^+_i\cap G^+_j$ and a vertex
$y\in G^-_i\cap G^-_j.$ Indeed, if such an edge $xy$ exists, since $i,j\in X$ and $i,j\notin Y$, we will obtain that $X=Y\cup \{ i,j\}$, contrary to the assumption that
the sets $X=\varphi(x)$ and $Y=\varphi(y)$ have the same size $m$. This shows that  $G^+_i\cap G^+_j$ and $G^-_i\cap G^-_j$ define different connected components of
$G(W_=(u,v))$. $\Box$

\medskip
Set $W'_=(u,v):=G^+_i\cap G^+_j$ and $W''_=(u,v):=G^-_i\cap G^-_j$ (we will use these notations even if some of  $G^+_i\cap G^+_j$ and $G^-_i\cap G^-_j$
are empty). By (2.7), they define the connected components of $G(W_=(u,v))$. Moreover, $W'_=(u,v)$ and $W''_=(u,v)$ are convex as the intersection of
two convex subsets of $G$.

\medskip
(2.8) {\it For an edge $e=(u,v)$ of $G$,  each of ${\mathcal W}'_{uv}$ and ${\mathcal W}''_{uv}$ constitutes complementary convex halfspaces of $G$.}

\medskip
{\it Proof.} Since $W(u,v)=G^+_i\cap G^-_j, W(v,u)=G^-_i\cap G^+_j$ by (2.6)  and $W'_=(u,v):=G^+_i\cap G^+_j$ and $W''_=(u,v)=G^-_i\cap G^-_j$ by their definition and (2.7),
we conclude that $W(u,v)\cup W'_=(u,v)=G^+_i\cap (G^-_j\cup G^+_j)=G^+_i$ and $W(v,u)\cup W''_=(u,v)=G^-_i\cap (G^+_j\cup G^-_j)=G^-_i$, hence they are complementary convex halfspaces of $G$.
Analogously one can show that $W(u,v)\cup W''_=(u,v)$ and  $W(v,u)\cup W'_=(u,v)$ are the complementary convex halfspaces $G^-_j$ and $G^+_j$ of $G$. $\Box$

\medskip
Summarizing, from (2.6)--(2.8) we obtain that a distance--preserving subgraph of a Johnson graph satisfies condition (WC).
Next we will show that $G$ satisfies condition (AGC). Let $b$ be a basepoint of $G$ and  $B:=\varphi(b)$.
Let ${\mathcal E}=\{ E_{\lambda}: \lambda\in \Upsilon\}$ be the set of  equivalence classes of $\theta_1$ with respect to $b$ and
let $\Sigma$  be the atom graph of $G$.

\medskip
(2.9) {\it For two vertical edges $e=uv$ and $e'=u'v'$ of $G$, $e\theta_1 e'$ iff $\gamma(e')=\gamma(e)$.}

\medskip
{\it Proof.} Let $\gamma(e)=\{ i,j\}$ with $V=U-i+j$. First suppose that $e\theta_1 e'$ and  let $u'\in W(u,v)$
and $v'\in W(v,u)$. By (2.8) and (2.2), $W(u,v)=W(u',v')$ and $W(v,u)=W(v',u')$. By (2.6), $W(u,v)=G^+_i\cap G^-_j$
and $W(v,u)=G^-_i\cap G^+_j$, hence $\{ i,j\}=U'\triangle V'=\gamma(u',v')$.  Conversely, if $\gamma(e')=\gamma(e)=\{ i,j\}$ and
$V'=U'-i+j,$ then $W(u',v')=G^+_i\cap G^-_j$ and $W(v',u')=G^-_i\cap G^+_j$, whence $W(u',v')=W(u,v), W(v',u')=W(v,u)$,
and by (2.2) we conclude that $e\theta_1 e'$. $\Box$ 

\medskip
 By (2.9), we can set $\gamma(E_{\lambda})=\{ i,j\}$ iff
$\gamma(e)=\{ i,j\}$ for an arbitrary (and thus for any) edge $e=uv\in E_{\lambda}$ with $b\in W(u,v)$. If $V=U-i+j,$ then set $\gamma^*(E_{\lambda}):=\{ j\}$.
Set $A:=\bigcup_{\lambda\in \Upsilon} \gamma^*(E_{\lambda})$. Let $D:=(A\cup B,F)$ be the graph having $A\cup B$ as the vertex--set and $i\in B$ and $j\in A$ are adjacent in $D$
(i.e., $ij\in F$) iff there exists  $E_{\lambda}\in {\mathcal E}$ such that $\gamma(E_{\lambda})=\{ i,j\}$. 

\medskip
(2.10) {\it $D=(A\cup B,F)$ is a bipartite graph with $B$ finite and the atom graph $\Sigma$ of $G$ is isomorphic to the line--graph $L(D)$ of $D$.}

\medskip
{\it Proof.} First we prove that $D$ is bipartite. Pick any vertex $j\in A$. By the definition of $A$ and (2.9), there exists an equivalence class $E_{\lambda}$ and a vertical edge $e=uv\in E_{\lambda}$ such that $\gamma(E_{\lambda})=\gamma(e)=\{ i,j\}$ and $\gamma^*(E_{\lambda})=\{ j\}$. If we assume without loss of generality that $b\in W(u,v)$, then this implies that $V=U-i+j.$
Now, if also $j\in B$, then this would imply that $b,v\in G^+_j$ and $u\in G^-_j$. Since $u\in I(v,b)$, this contradicts the convexity of the set $G^+_j$. Hence $A\cap B=\emptyset$. Since by definition of $D$, any edge of $D$ is running between a vertex of $B$ and a vertex of $A$, the graph $D$ is bipartite. Finally, since $B=\varphi(b),$ the set $B$ is finite.

Now we will prove that the graphs  $\Sigma$ and $L(D)$ are isomorphic, namely, that $\gamma$ is an isomorphism between  $\Sigma$ and $L(D)$. First, consider any edge $E_{\lambda}E_{\lambda'}$ of $\Sigma$ and pick arbitrary edges $e=uv\in E_{\lambda}, e'=u'v'\in E_{\lambda'}$, where  $b\in W(u,v)\cap W(u',v')$. We will show that $|\gamma(e)\cap \gamma(e')|=1$. Suppose that $\gamma(e)=\gamma(E_{\lambda})=\{ i,j\},$ where $i\in B$ and $j\in A$. Then $V=U-i+j$. By (2.9), $\gamma(e)\ne \gamma(e')$.  By (2.3), the edge $u'v'$ has one end in $W(u,v)\cup W(v,u)$ and another end in $W_=(u,v)$. We distinguish two cases:

\medskip\noindent
{\bf Case 1:} $u'\in W_=(u,v)$. By (2.3),  necessarily $v'\in W(v,u).$  Hence $v'\in G^+_j$. If $u'\in G^-_j,$ then $U'\triangle V'=\{ i',j\}$ for some $i'\in B$ and we conclude that $\gamma(E_{\lambda'})=\gamma(e')=\{ i',j\}$, thus $E_{\lambda}E_{\lambda'}$  corresponds to  edges $ij$ and $i'j$ of $D$. Hence, we can suppose that $u'\in G^+_j$. Since $u'\in W_=(u,v),$ by (2.7) we conclude that  $u'\in G^+_i\cap G^+_j$. On the other hand,
since $v\in I(u,v')$ and $u\in G^+_i,v\in G^-_i,$ from the convexity of $G^+_i$ we deduce that $v'\in G^-_i$. Hence $V'=U'-i+j'$ and $\gamma(E_{\lambda'})=\gamma(e')=\{ i,j',\}$, thus  $E_{\lambda}E_{\lambda'}$ corresponds to edges $ij$ and $ij'$ of $D$.

\medskip\noindent
{\bf Case 2:} $v'\in W_=(u,v)$. Then $u'\in W(u,v)$ by (2.3).  By (2.6), $u'\in G^+_i\cap G^-_j$. If $v'\notin G^+_i$, then $V'=U'-i+j'$, and consequently $\gamma(E_{\lambda'})=\gamma(e')=\{ i,j',\}$, whence   $E_{\lambda}E_{\lambda'}$  corresponds to  edges $ij$ and $ij'$ of $D$. On the other hand, if $v'\in G^+_i$, since $v'\in W_=(u,v)$ by (2.7) we conclude that
$v'\in G^+_i\cap G^+_j$. On the other hand, since $u\in I(v,u')$ and $v\in G^+_j,u\in G^-_j,$ from the convexity of $G^+_j$ we deduce that  $u'\in G^-_j$. Hence $V'=U'-i'+j$ and $\gamma(E_{\lambda'})=\gamma(e')=\{ i',j,\}$, thus  $E_{\lambda}E_{\lambda'}$  corresponds to  edges $ij$ and $i'j$ of $D$.

Conversely, pick an edge of $L(D)$ corresponding to two incident edges $ij$ and $i'j'$ of $D$, where $i,i'\in B$ and $j,j'\in A$. By definition of  $D$, there exist two equivalence
classes $E_{\lambda}$ and $E_{\lambda'}$ of ${\mathcal E}$ such that $\gamma(E_{\lambda})=\{ i,j\}$ and  $\gamma(E_{\lambda'})=\{ i',j'\}.$ We will prove that $E_{\lambda}$ and $E_{\lambda'}$
are adjacent in $\Sigma$. Pick two edges $e=uv\in E_{\lambda}$ and $e'=u'v'\in E_{\lambda'}$
such that $b\in W(u,v)\cap W(u',v')$.  Then $V=U-i+j$ and $V'=U'-i'+j'$. Again, we distinguishes two cases:

\medskip\noindent
{\bf Case 3:} $i=i'$ and $j\ne j'$. The convexity of $G^+_i$ implies that $u'\notin W(v,u)$. Since $v,v'\in G^-_i$ and $u\in G^+_i$, the convexity of $G^-_i$ implies that $v'\notin W(u,v)$. First suppose that  $u'\in W(u,v)$. If $v'\in W(v,u),$ then $e\theta_1e'$ and from (2.9) we deduce that $\gamma(e)=\gamma(e')$, thus $j=j'$, which is impossible. Hence $v'\in W_=(u,v)$ and
by (2.3) we have $e\theta_2 e'$. By (2.4),  $E_{\lambda}E_{\lambda'}$ is an edge of $\Sigma$. Now suppose that $u'\in W_=(u,v)$. If $v'\in W(v,u),$ then analogously  we can conclude that $E_{\lambda}$ and $E_{\lambda'}$ are adjacent in $\Sigma$. Therefore it remains to consider the case where $v'\in W_=(u,v)$. Since $u'$ and $v'$ are adjacent in $G$, they belong to the same connected component of $W_=(u,v)$. Since $u'\in G^+_i$ and $v'\in G^-_i$, this contradicts (2.7).

\medskip\noindent
{\bf Case 4:} $j=j'$ and $i\ne i'$. The convexity of $G^+_j$ implies that $v'\notin W(u,v)$. Since $u,u'\in G^-_j$ and $v\in G^+_j$, the convexity of $G^-_j$ implies that $u'\notin W(v,u)$. First suppose that $u'\in W(u,v)$. Then either $v'\in W(v,u)$ and $e\theta_1 e',$ which is impossible, or $v'\in W_=(u,v)$ and we can proceed  as in Case 3. Finally suppose that $u'\in W_=(u,v)$. Again, if $v'\in W(v,u)$ we can conclude the proof as in Case 3. On the other hand, if $v'\in W_=(u,v),$ then $u'$ and $v'$ belong to the same connected component of $W_=(u,v)$. Since $u'\in G^-_j$ and $v'\in G^+_j$, this contradicts (2.7). This establishes that the atom graph $\Sigma$ is the line--graph of the bipartite graph $D=(A,B;F)$. $\Box$

\subsection{Sufficiency}
We will prove the converse implication that a graph $G=(V,E)$ satisfying the conditions (WC) and (AGC) is isometrically embeddable into a Johnson graph. By (WC) and (2.2), the relation $\theta_1$ is transitive, thus the atom graph $\Sigma$ in condition (AGC) is well--defined. Let ${\mathcal E}=\{ E_{\lambda}: \lambda\in \Upsilon\}$ denote the vertex--set of $\Sigma$, i.e., the set of equivalence classes of $\theta_1$. By (2.4), $E_{\lambda}$ and $E_{\lambda'}$ are adjacent in $\Sigma$ iff $e\theta_2 e'$ for arbitrary edges $e\in E_{\lambda}$ and $e'\in E_{\lambda'}$. Suppose that $\Sigma$ is the line--graph of a bipartite graph $D=(A,B;\Upsilon)$, where the set $B$ is finite, say $|B|=m$. Set $\Lambda:=A\cup B$.  Let $b$ be an arbitrary vertex of $G$ and let $T$ be a spanning tree of $G$ rooted at $b$ such that $d_T(b,v)=d_G(b,v)$ for all $v\in V$. For a vertex $v$ of $G$, denote by $P_v$ the path of $T$ connecting $b$ to $v$.
Each edge $e=uv$ of $T$ is a vertical edge of $G$: if $u$ is the father of $v$ in $T$, then $d_G(b,v)=d_G(b,u)+1$. Therefore  $e$  belongs to some $E_{\lambda}\in {\mathcal E}$. Then $E_{\lambda}$ is a vertex of $\Sigma$ corresponding to some edge $ij$ of $D$ with $i\in B$ and $j\in A$. Set $\gamma(e)=\gamma(E_{\lambda}):=\{ i,j\}$ and call it the  {\it label} of the edge $e$.
We construct an isometric embedding $\varphi$ of $G$ into $J(m,\Lambda)$ inductively, according to the distance $d_G(v,b)$ from the current vertex $v$ to $b$.  Set $\varphi(b):=B$. In the assumption that $uv$ is an edge of $T$ with $d_G(v,b)=d_G(u,b)+1$ such that $\varphi(u)$  has been defined (let $U:=\varphi(u)$), we define $\varphi(v)$ by setting $V=\varphi(v):=U-i+j$.

\medskip
(2.11) {\it For any two vertical edges $e=uv$ and $e'=u'v'$ of $G$ with $b\in W(u,v)\cap W(u',v'),$  $|\gamma(e)\cap \gamma(e')|=<e,e'>$.}

\medskip
{\it Proof.} First note that $\gamma(e)=\gamma(e')$ iff  $e$ and $e'$ belong to the same equivalence class of $\mathcal E$, i.e.,  $<e',e''>=2$ by (2.2). If  $|\gamma(e)\cap \gamma(e')|=1$, then
$e$ and $e'$ belong to different classes $E_{\lambda},E_{\lambda'}$ of $\mathcal E$ corresponding to two adjacent vertices of $\Sigma$, whence $<e,e'>=1$ by (2.3). Conversely, let $<e,e'>=1$. Then $e$ and $e'$ belong to different equivalent classes $E_{\lambda}$ and $E_{\lambda'}$ of $\mathcal E$. By (2.4),  $E_{\lambda}$ and $E_{\lambda'}$ are adjacent vertices of $\Sigma$. Therefore $E_{\lambda'}$ and $E_{\lambda''}$ correspond to two incident edges of  $D$, hence $|\gamma(e)\cap \gamma(e')|=|\gamma(E_{\lambda})\cap \gamma(E_{\lambda'})|=1$. Consequently,
$|\gamma(e)\cap \gamma(e')|=0$ iff the edges $e'$ and $e''$ are not in relation $\theta_1$ or $\theta_2$, i.e., $<e,e'>=0$.  $\Box$

\medskip
(2.12) {\it For any vertex $v$ of $G$,  the labels of all edges of the path $P_v$ are pairwise disjoint, $|V\triangle B|=2d_G(v,b)$, and $|V|=m$.}

\medskip
{\it Proof.} Pick any two edges $e$ and $e'$ of $P_v$. Since  $P_v$ is a shortest path of $G$,
$e$ and $e'$ are not in relation $\theta_1$ or $\theta_2$.  By (2.11), $|\gamma(e)\cap \gamma(e')|=<e,e'>=0$.  Let $u$ be the neighbor of $v$ in $P_v$.
Since $\gamma(uv)=\{ i,j\}$ and $i\in B$, from the definition of $\varphi$ we conclude that $i$ and $j$ do not occur in the labels of the edges of the path $P_u$.
Hence $i\in U$ and $j\notin U$. Since $V=U-i+j$, we obtain $|V|=|U|=m$. While moving along  $P_v$ from $b$ to $v$, when we traverse an edge we remove an
element of $B$ and add an element of $A$. Since the removed and added elements are pairwise distinct,  $|V\triangle B|$ is twice the length of $P_v$, i.e., $2d_G(v,b)$. $\Box$

\medskip
Hence $\varphi$ maps the vertices of $G$ to vertices of  $J(m,\Lambda)$. It remains to show that $\varphi$ is an isometric embedding, i.e., for any two vertices $x,y$ of $G$, $|\varphi(x)\triangle \varphi(y)|=|X\triangle Y|=2d_G(x,y)$. For technical conveniences,  we proceed as follows. First, we define the mapping $\psi$ from the hypercube $H(\Lambda)$ to itself by setting  $X'=\psi(X):=X\triangle B$ for each finite set $X$.
Since $X\triangle Y=(X\triangle B)\triangle (Y\triangle B)=X'\triangle Y'$, $\psi$ is an isometry of $H(\Lambda)$. Since $(X\triangle B)\triangle B=X$, $\psi$ is an involution.
Note that $\psi$ maps isometrically the half--cube containing the vertex--set of $J(m,\Lambda)$ to the half--cube containing $B'=\emptyset$. Therefore, for two vertices $x,y$ of $G$,
the equality $|X\triangle Y|=2d_G(x,y)$ holds iff $|X'\triangle Y'|=2d_G(x,y)$. We will prove this second equality following quite literally an analogous proof from
\cite{Shp} and \cite[Subsection 21.3(v)]{DeLa}.

\medskip
(2.13) {\it For any two vertices $x,y$ of $G$, $2d_G(x,y)=|X'\triangle Y'|=|X\triangle Y|$. Consequently, $\varphi$ is an isometric embedding of $G$ into the Johnson graph $J(m,\Lambda)$.}

\medskip
{\it Proof.} First, summarizing (2.11) and (2.12), we conclude that for any vertex $v$ of $G$, we have $|V'|=|V'\triangle B'|=|V\triangle B|=2d_G(v,b)$ and that $V'$ is the union of  labels $\gamma(e)$ of all edges on the path $P_v$ of $T$. Therefore, we can suppose that $x$ and $y$ are different from $b$. We prove  $2d_G(x,y)=|X'\triangle Y'|$ by induction on the sum $d_G(b,x)+d_G(b,y)$. Let $x_1$ be the neighbor of $x$ in $P_x$ and $y_1$ be the neighbor of $y$ in $P_y$. Set $e:=x_1x$ and $e':=y_1y$. Let  $X',Y',X'_1,$ and $Y'_1$ be the images of $X,Y,X_1,$ and $Y_1$ (and thus of $x,y,x_1$ and $y_1$) under the map $\psi$.  By the induction assumption, we have
\begin{align*}
|X'_1\triangle Y'_1|=2d_G(x_1,y_1), |X'_1\triangle Y'|=2d_G(x_1,y), |X'\triangle Y'_1|=2d_G(x,y_1). \tag{4}
\end{align*}
Since $2|V'\cap W'|=|V'|+|W'|-|V'\triangle W'|$ and $V'\triangle B'=V', W'\triangle B'=W'$ hold for arbitrary finite subsets $V'$ and $W'$ of $\Lambda$, from (4) we obtain:
\begin{align*}\begin{cases}
|X'_1\cap  Y'_1|=d_G(b,x_1)+d_G(b,y_1)-d_G(x_1,y_1),\\
|X'_1\cap  Y'|=d_G(b,x_1)+d_G(b,y)-d_G(x_1,y),\\
|X'\cap  Y'_1|=d_G(b,x)+d_G(b,y_1)-d_G(x,y_1). \tag{5}
\end{cases}
\end{align*}
By (2.12)  and the definition of $\psi$, we have $X'=X'_1\cup \gamma(e), Y'=Y'_1\cup \gamma(e')$, thus
\begin{align*}
|X'\cap Y'|=|X'_1\cap Y'|+|Y'_1\cap X'|-|X'_1\cap Y'_1|+|\gamma(e)\cap \gamma(e')|. \tag{6}
\end{align*}
On the other hand,
\begin{align*}
|X'\triangle Y'|=|X'|+|Y'|-2|X'\cap Y'|=2d_G(b,x)+2d_G(b,y)-2|X'\cap Y'|.\tag{7}
\end{align*}
Substituting  the three equalities of (5) in (6) and then inserting the result in (7), we obtain
\begin{align*}
|X'\triangle Y'|=2(d_G(x,y_1)+d_G(x_1,y)-d_G(x_1,y_1)-|\gamma(e)\cap \gamma(e')|).\tag{8}
\end{align*}
By (2.11) and (3), $|\gamma(e)\cap \gamma(e')|=<e,e'>=d_G(x,y_1)+d_G(x_1,y)-d_G(x,y)-d_G(x_1,y_1)$.
Therefore, the right--hand side of (8) equals to
$2d_G(x,y),$ establishing the required equality $|X'\triangle Y'|=2d_G(x,y).$ $\Box$

\medskip
This concludes the proof of Theorem 1.

\subsection{Proof of Corollary 1} In view of the results of \cite{Mau} and \cite{ChChOs}, it suffices to show that
the wallspace condition (WC) implies the positioning condition (PC). Suppose by way of contradiction that $G$ satisfies
(WC) but contains a vertex $b$ and a square $u_1u_2u_3u_4$ such that  $d_G(b,u_1)+d_G(b,u_3)<d_G(b,u_2)+d_G(b,u_4).$
Suppose that $k:=d_G(b,u_1)\le d_G(b,u_3)$
and $d_G(b,u_2)\le d_G(b,u_4)$. Then obviously $d_G(b,u_4)=k+1.$ We have three possibilities: (a) $d_G(b,u_2)=d_G(b,u_3)=k$,
(b) $d_G(b,u_2)=d_G(b,u_3)=k+1,$ and (c) $d_G(b,u_2)=k+1, d_G(b,u_3)=k$.

In case (a), we have $b\in W_=(u_1,u_2), u_4\in W(u_1,u_2),$
and $u_3\in W(u_2,u_2)$. Since $u_3\in I(b,u_4)$, this contradicts the convexity of $W(u_1,u_2)\cup W'_=(u_1,u_2),$
where $W'_=(u_1,u_2)$ is the connected component of $W_=(u_1,u_2)$ containing $b$. In case (b), we have $b\in W_=(u_2,u_3), u_1\in W(u_2,u_3)$,
and $u_4\in W(u_3,u_2)$. Since $u_1\in I(b,u_4)$, this contradicts the convexity of   $W(u_3,u_2)\cup W'_=(u_2,u_3),$ where
$W'_=(u_2,u_3)$ is the connected component of $W_=(u_2,u_3)$ containing $b$. Finally, in case (c), we have $b,u_4\in W(u_1,u_2)$ and $u_3\in W(u_2,u_1)$.
Since $u_3\in I(b,u_4),$ this contradicts the convexity of $W(u_1,u_2)$. This concludes the proof of Corollary 1.

\section{Concluding remarks}

\medskip
(3.1) Condition (WC) and the assertions (2.7) and (2.8) completely describe the wallspace $\frak W$ defining an isometric embedding of a graph $G$ into a Johnson graph: $\frak W$ consists of all distinct walls of the form ${\mathcal W}'_{uv}=\{ W(u,v)\cup W'_=(u,v), W(v,u)\cup W''_=(u,v)\}$ and ${\mathcal W}''_{uv}=\{ W(u,v)\cup W''_=(u,v), W(v,u)\cup W'_=(u,v)\}$, where  $uv$ runs over the vertical edges of $G$.  The positive and negative signs of the halfspaces of the walls ${\mathcal W}'_{uv}, {\mathcal W}''_{uv}$  are determined by the choice of the basepoint $b$ and the encoding of the atom graph $\Sigma$ as the line--graph of a bipartite graph. Notice that if $W(u,v)\cup W'_=(u,v)$ is a negative halfspace and its complement  $W(v,u)\cup W''_=(u,v)$ is a positive halfspace of ${\mathcal W}'_{uv}$, then  $W(u,v)\cup W''_=(u,v)$ will be a positive halfspace  and $W(v,u)\cup W'_=(u,v)$ will be a negative halfspace of ${\mathcal W}''_{uv}$.

\medskip
(3.2) If $G$ is a distance--preserving subgraph of a hypercube, then  $G$  satisfies condition (WC): $W_=(u,v)=\emptyset$ for each edge $uv$, hence the two walls ${\mathcal W}'_{uv},{\mathcal W}''_{uv}$ separating $u$ and $v$ coincide. For any basepoint $b$, all edges are vertical and no edges are in relation $\theta_2$  because $G$ is bipartite. Therefore $\Sigma$ consists of isolated vertices, one for each equivalence class of $\mathcal E$, i.e., $\Sigma$ is the line--graph of a matching. Consequently, $G$ is isometrically embeddable into a Johnson graph iff $\Sigma$ is finite (which is equivalent to the finiteness of $G$).

\medskip
(3.3) It is instructive to apply Theorem 1 to the Petersen graph $P_{10}$. For any edge $uv$, $W_=(u,v)$ has two connected components, each consisting of an  edge.  The
two walls ${\mathcal W}'_{uv}, {\mathcal W}''_{uv}$ separating $u$ and $v$ consist of two disjoint 5--cycles,
which are convex subgraphs of $P_{10}$. For any basepoint $b$,  $\Sigma$ has 9 vertices and is 4--regular. It can be easily seen that $\Sigma$ is the line--graph of $K_{3,3}$. Assuming
that $K_{3,3}$ has the vertices $1,2,3$ in one part and $4,5,6$ in another part and $b$ is labeled by $123$, the labels of the remaining  vertices of $P_{10}$ are $125,136,234,246,345,356,256,146,$ and $145$.

\begin{figure}[h]
\begin{center}
\includegraphics{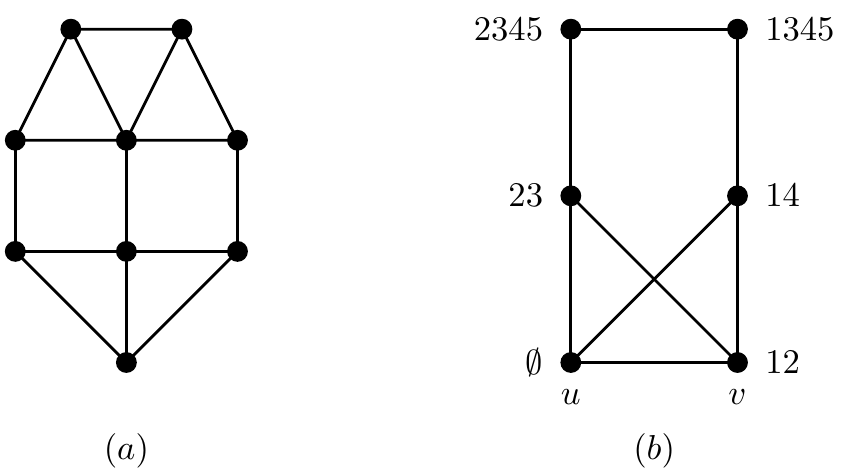}%
\end{center}
\caption{Two graphs $G_1$ and $G_2$ not isometrically embeddable into a Johnson graph.}
\label{figure1}
\end{figure}

\medskip
(3.4) One can ask if condition (AGC)  can be replaced by the link condition (LC): the
neighborhoods of all vertices of $G$ induce  line--graphs of bipartite graphs. (LC) is necessary for an isometric embedding  into a Johnson graph. However, the graph
$G_1$ in Figure \ref{figure1}(a) shows that (LC)
cannot replace (AGC). Indeed, $G_1$ satisfies (WC) and  (LC). On the other hand, if we pick the bottom vertex of $G_1$
as the basepoint $b$, then the atom graph of $G_1$ is not the line--graph of a bipartite graph.

\medskip
(3.5)  Figure \ref{figure1}(b) presents a graph $G_2$ and its isometric embedding into the half--cube $\frac{1}{2}H(\Lambda)$ with $\Lambda=\{ 1,2,3,4,5\}$. On the other hand,
$G_2$ is not a distance--preserving subgraph of a Johnson graph: for edge $e=uv$ as in the figure, $W_=(u,v)$ has a single connected component consisting of all vertices
of $G_2$ except $u$ and $v$. However, $W_=(u,v)$ and its unions with $W(u,v)$ and $W(v,u)$ are not convex, thus $G_2$ violates condition (WC).
This shows that  in the  wallspaces defining isometric embeddings of graphs into half--cubes the walls can split the connected components of $W_=(u,v)$, thus leaving
less chances for a structural characterization of distance--preserving subgraphs of half--cubes (recall \cite[Problem 21.4.1]{DeLa}).

\medskip
(3.6) Many important graph classes are isometrically embeddable into hypercubes. Among them we can mention the graphs of regions of hyperplane arrangements in ${\mathbb R}^d$,
and, more generally, the tope graphs of oriented matroids \cite{BjLVStWhZi} and of complexes of oriented matroids \cite{BaChKn}, the median graphs (alias 1--skeleta of CAT(0) cube complexes) \cite{BaCh_survey,ImKl}, the graphs of lopsided sets
\cite{BaChDrKo1,La}, and the 1--skeleta of CAT(0) Coxeter zonotopal complexes \cite{HaPa} (for examples of $\ell_1$--graphs see the book \cite{DeLa} and the survey \cite{BaCh_survey}).

\medskip\noindent
{\bf Acknowledgements.} I am indepted to the anonymous referees for careful reading of the first version and several corrections.

\bibliographystyle{amsalpha}

\end{document}